\documentclass[12pt]{article} 

\usepackage{amsmath}
\usepackage{amssymb} 
\usepackage{amscd} 
\usepackage{theorem} 
\usepackage{a4} 

\def\oplusinf{\mathop{\oplus}} 

\def\im{{\mbox{Im}}}

\def\ker{{\mbox{Ker}}}

\newtheorem{theorem}{THEOREM}
\newtheorem{lemma}{LEMMA}

\begin{document}

\baselineskip=0.7cm 
\begin{center} 
{\Large\bf GENERALIZED COHOMOLOGY}

\vspace{3mm}

{\Large\bf FOR IRREDUCIBLE TENSOR FIELDS}

\vspace{3mm}

{\Large\bf  OF MIXED YOUNG SYMMETRY TYPE} 

\end{center} 
\vspace{0.75cm}

\begin{center} Michel DUBOIS-VIOLETTE \\
\vspace{0.3cm} {\small Laboratoire de Physique Th\'eorique
\footnote{Unit\'e Mixte de Recherche du Centre National de la
Recherche Scientifique - UMR 8627}\\ Universit\'e Paris XI,
B\^atiment 210\\ F-91 405 Orsay Cedex, France\\
patricia$@$lyre.th.u-psud.fr}\\ and\\
Marc HENNEAUX\\
\vspace{0.3cm} {\small Physique Th\'eorique et Math\'ematique\\ 
Universit\'e Libre de Bruxelles\\
Campus Plaine C.P. 231\\ B-1050 Bruxelles, Belgique
\footnote{Also at Centro de Estudios Cient\'{\i}ficos de Santiago,
Casilla 16443, Santiago 9, Chile}
\\
henneaux$@$ulb.ac.be}

\end{center} \vspace{1cm}

\begin{center} \today \end{center}

\vspace {1cm}

\noindent LPT-ORSAY 99/56\\
ULB-TH 99/12

\newpage

\begin{abstract}
We construct $N$-complexes of non completely antisymmetric
irreducible tensor fields on $\mathbb R^D$ generalizing thereby
the usual complex $(N=2)$ of differential forms. 
These complexes arise naturally in the description
of higher spin gauge fields.  Although, for  
$N\geq 3$, the generalized cohomology of these $N$-complexes is
non trivial, we prove a generalization of the Poincar\'e lemma.
Several results which appeared in various contexts are shown
to be particular cases of this generalized Poincar\'e lemma. 
\end{abstract}
\newpage

\section{Introduction}

Our aim in this letter is to set up differential tools for
irreducible tensor fields on $\mathbb R^D$ which generalize the
calculus of differential forms. By an irreducible tensor field on
$\mathbb R^D$, we here mean, a smooth mapping $x\mapsto T(x)$ of
$\mathbb R^D$ into a vector space of (covariant) tensors of given
Young symmetry. We recall that this implies that the
representation of $GL_D$ in the corresponding space of tensors is
irreducible.

We first introduce a generalization of 
the familiar exterior derivative that
satisfies, instead of $d^2 = 0$,
the nilpotency condition $d^N = 0$ for some integer 
$N \geq 2$  that depends on the Young symmetry type
of the tensor fields under consideration.
We then analyse the
generalized (co)homologies $H^{(k)} \equiv Ker d^k /Im d^{N-k}$,
$(k=1, \cdots, N-1)$ of these nilpotent endomorphisms in the sense of
\cite{D-V}, \cite{D-V2}, \cite{D-V3}, \cite{D-VK}, \cite{Kap}, \cite{KW} and establish an analog of
the Poincar\'e lemma.  

The nilpotent endomorphisms introduced here have various physical 
applications.
They naturally arise, for instance, in the theory of 
higher spin gauge fields.  They also encompass 
conservation laws involving symmetric tensors.
This is discussed at the end of the letter.

An expanded version of this letter, with further developments and
detailed proofs, will appear elsewhere \cite{D-VH}.

\section{Definitions}

Throughout the following $(x^\mu)=(x^1,\dots,x^D)$ denotes the
canonical coordinates of $\mathbb R^D$ and $\partial_\mu$ are the
corresponding partial derivatives which we identify with the
corresponding covariant derivatives associated to the canonical
flat linear connection of $\mathbb R^D$. Thus, for instance, if $T$
is a covariant tensor field of degree $p$ on $\mathbb R^D$ with
components $T_{\mu_1\dots\mu_p}(x)$, then $\partial T$ denotes
the covariant tensor field of degree $p+1$ with components
$\partial_{\mu_1}T_{\mu_2\dots\mu_{p+1}}(x)$. The operator
$\partial$ is a first-order differential operator which increases
by one the tensorial degree.

In this context, the space $\Omega(\mathbb R^D)$ of differential
forms on $\mathbb R^D$ is the graded vector space of (covariant)
antisymmetric tensor fields on $\mathbb R^D$ with graduation
induced by the tensorial degree whereas the exterior differential
$d$ is the composition of the above $\partial$ with
antisymmetrisation, i.e.
\begin{equation}
d={\mathbf A}_{p+1}\circ \partial : \Omega^p(\mathbb R^D)\rightarrow
\Omega^{p+1}(\mathbb R^D)
\label{eq1}
\end{equation}
where ${\mathbf A}_p$ denotes the antisymmetrizer on tensors of degree $p$.
The Poincar\'e lemma asserts that the cohomology of the complex
$(\Omega(\mathbb R^D),d)$ is trivial, i.e. that one has
$H^p(\Omega(\mathbb R^D))=\ker(d:\Omega^p(\mathbb R^D)\rightarrow
\Omega^{p+1}(\mathbb R^D))/d(\Omega^{p-1}(\mathbb R^D))=0$,
$\forall p\geq 1$ and $H^0(\Omega(\mathbb
R^D))=\ker(d:\Omega^0(\mathbb R^D)\rightarrow \Omega^1(\mathbb
R^D))=\mathbb R$.

{}From the point of view of Young symmetry, antisymmetric
tensors correspond to Young diagrams (partitions) described 
by one column of cells, i.e. the space of values of
$p$-forms corresponds to one column of $p$ cells, $(1^p)$,
whereas ${\mathbf A}_p$ is the associated Young symmetrizer.

There is a relatively easy way to generalize the pair
$(\Omega(\mathbb R^D),d)$ which we now describe. Let
$Y=(Y_p)_{p\in \mathbb N}$ be a sequence of Young diagrams
such that the number of cells of $Y_p$ is $p$, $\forall p\in
\mathbb N$ (i.e. such that $Y_p$ is a partition of the integer
$p$ for any $p$).
We define $\Omega^p_Y(\mathbb R^D)$ to be
the vector space of smooth covariant tensor fields of degree
$p$ on $\mathbb R^D$ which have the Young symmetry type $Y_p$
and we let $\Omega_Y(\mathbb R^D)$ be the graded vector space
$\displaystyle{\oplusinf_p}\Omega^p_Y(\mathbb R^D)$.  We then
generalize the exterior differential by setting $d={\mathbf
Y}\circ \partial$, i.e.
\begin{equation}
d={\mathbf Y}_{p+1}\circ\partial :\Omega^p_Y(\mathbb
R^D)\rightarrow \Omega^{p+1}_Y(\mathbb R^D)
\end{equation}
where ${\mathbf Y}_p$ is now the Young symmetrizer on tensor
of degree $p$ associated to the Young symmetry $Y_p$. This $d$
is again a first order differential operator which is of
degree one, (i.e. it increases the tensorial degree by one),
but now, $d^2\not= 0$ in general. Instead, one has the
following result.

\begin{lemma}
Let $N$ be an integer with $N\geq 2$ and assume that $Y$ is
such that the number of columns of the Young diagram $Y_p$ is
strictly smaller than $N$ (i.e. $\leq N-1$) for any $p\in
\mathbb N$. Then one has $d^N=0$.
\end{lemma}

In fact the indices in one column are antisymmetrized and $d^N\omega$ involves necessarily at least two
partial derivatives $\partial$ in one of the columns since
there are $N$ partial derivatives involved and at most $N-1$
columns.

Thus if $Y$ satisfies the condition of Lemma 1,
$(\Omega_Y(\mathbb R^D),d)$ is a $N$-{\sl complex} (of cochains)
\cite{Kap}, \cite{D-V}, \cite{D-VK}, \cite{KW}, \cite{D-V2},
i.e. here a graded vector space equipped with an endomorphism
$d$ of degree 1, its $N$-{\sl differential}, satisfying $d^N=0$. Concerning $N$-complexes, we
shall use here the notations and the results \cite{D-V2}.

Notice that $\Omega^p_Y(\mathbb R^D)=0$ if the first column of
$Y_p$ contains more than $D$ cells and that therefore, if $Y$
satisfies the condition of Lemma 1, then $\Omega^p_Y(\mathbb
R^D)=0$ for $p>(N-1)D$. 

One can also define a graded bilinear product on
$\Omega_Y(\mathbb R^D)$ by setting 
\begin{equation}
(\alpha\beta)(x)={\mathbf Y}_{a+b}(\alpha(x)\otimes \beta(x))
\label{eq3}
\end{equation}
for $\alpha\in \Omega^a_Y(\mathbb R^D)$, $\beta\in
\Omega^b_Y(\mathbb R^D)$ and $x\in \mathbb R^D$. This product is
by construction bilinear with respect to the $C^\infty(\mathbb
R^D)$-module structure of $\Omega_Y(\mathbb R^D)$ (i.e. with
respect to multiplication by smooth functions). It is worth
noticing here that one always has $\Omega^0_Y(\mathbb
R^D)=C^\infty(\mathbb R^D)$.

\section{The $N$-complexes $(\Omega_N(\mathbb R^D), d)$}

In this letter, we shall not stay at this level of generality but,
for each $N\geq 2$ we shall choose a maximal $Y$, denoted by
$Y^N=(Y^N_p)_{p\in\mathbb N}$, satisfying the condition of lemma
1. The Young diagram with $p$ cells $Y^N_p$ is defined in the
following manner: write the division of $p$ by $N-1$, i.e. write
$p=(N-1)n_p+r_p$ where $n_p$ and $r_p$ are (the unique) integers
with $0\leq n_p$ and $0\leq r_p\leq N-2$ ($n_p$ is the quotient
whereas $r_p$ is the remainder), and let $Y^N_p$ be the Young
diagram with $n_p$ rows of $N-1$ cells and the last row with
$r_p$ cells (if $r_p\not= 0$). One has
$Y^N_p=((N-1)^{n_p},r_p)$, that is we fill the rows maximally.
We shall denote $\Omega_{Y^N}(\mathbb R^D)$ and
$\Omega^p_{Y^N}(\mathbb R^D)$ by $\Omega_N(\mathbb R^D)$ and
$\Omega^p_N(\mathbb R^D)$. It is clear that $(\Omega_2(\mathbb
R^D),d)$ is the usual complex $(\Omega(\mathbb R^D),d)$ of
differential forms on $\mathbb R^D$. The $N$-complex
$(\Omega_N(\mathbb R^D),d)$ will be simply denoted by
$\Omega_N(\mathbb R^D)$. 

We call the Young diagrams $Y^N_p$
with $p=(N-1)n_p$ ``well-filled diagrams".
These are rectangular diagrams with $n_p$ rows of
$N-1$ cells each.

We recall \cite{D-V2} that the
(generalized) cohomology of the $N$-complex $\Omega_N(\mathbb
R^D)$ is the family of graded vector spaces $H_{(k)}(\Omega_N(\mathbb R^D))$ $k\in
\{1,\dots,N-1\}$  defined by
$H_{(k)}(\Omega_N(\mathbb R^D))=\ker(d^k)/\im(d^{N-k})$, i.e.
$H_{(k)}(\Omega_N(\mathbb
R^D))=\displaystyle{\oplusinf_p}H^p_{(k)}(\Omega_N(\mathbb
R^D))$ with 
\[
H^p_{(k)}(\Omega_N(\mathbb
R^D))=\ker(d^k:\Omega^p_N(\mathbb R^D)\rightarrow \Omega^{p+k}_N
(\mathbb R^D))/d^{N-k}(\Omega^{p+k-N}(\mathbb R^D)).
\]

It is easy to write down explicit formulas in terms of components.
Consider for instance the case $N = 3$, for which the
relevant Young diagrams are those with two colums, one
of length $k$ and the second of length $k-1$ or $k$.  A tensor field in
$\Omega_3(\mathbb R^D)$ is a scalar $T$ in tensor degree $0$,
a vector $T_\alpha$ in tensor degree $1$, a symmetric tensor
$T_{\alpha \beta}$ in tensor degree $2$.  In tensor degree
$2k-1$ ($k \geq 2$), it is described by
components $T_{\alpha_1 \dots \alpha_k \beta_1 \dots \beta_{k-1}}$
with the Young symmetry of the diagram with $k-1$ rows
of length $2$ and one row of length $1$, while
in even tensor degree
$2k$, it is described by components 
$T_{\alpha_1 \dots \alpha_k \beta_1 \dots \beta_k}$ 
with the Young symmetry of the well-filled rectangular diagram
with $k$ rows of length $2$.
The components of $dT$ are respectively proportional to
$\partial_\alpha T$, $\partial_{(\alpha} T_{\beta)}$,
$\partial_{[\alpha_1} T_{\alpha_2] \beta}$ and 
$ T_{\alpha_1 \dots \alpha_k [\beta_2 \dots \beta_{k}},_{\beta_1]}
+ T_{\beta_1 \dots \beta_k [\alpha_2 \dots \alpha_{k}},_{\alpha_1]}$
or $ \partial_{[\alpha_{1}} T_{\alpha_2 \dots \alpha_{k+1}]
\beta_1 \dots \beta_k}$, where the comma stands for the partial
derivative, $(\dots)$ for symmetrization and $[\dots]$ for antisymmetrization.  It is obvious that $d^3 =0$ since all terms in $d^3T$
involves one antisymmetrization over partial derivatives.

\section{Generalized Poincar\'e lemma}

The following statement is our generalization of the Poincar\'e
lemma.

\begin{theorem}
One has $H^{(N-1)n}_{(k)}(\Omega_N(\mathbb R^D))=0$, $\forall
n\geq 1$ and $H^0_{(k)}(\Omega_N(\mathbb R^D))$ is the space of
real polynomial functions on $\mathbb R^D$ of degree strictly
less than $k$ (i.e. $\leq k-1$) for $k\in\{1,\dots,N-1\}$.
\end{theorem}

This statement reduces to the Poincar\'e lemma for $N=2$ but it
is a nontrivial generalization for $N\geq 3$ in the sense that
the spaces $H^p_{(k)}(\Omega_N(\mathbb R^D))$
are nontrivial for $p\not=(N-1)n$ and, in fact, are generically
infinite dimensional for $D\geq 3$, $p\geq N$.

The second part of the theorem is obvious since the condition
$d^k f =0$ simply states that the derivatives of order $k$ of
$f$ all vanish (and there is no quotient to be taken since $f$
is in degree $0$).
The proof of the first part of the
theorem, which asserts that there is no cohomology
for well-filled diagrams,  proceeds by introducing an appropriate
generalized homotopy \cite{D-V2}.  By inner contraction with the
vector $x^\mu$, one can easily construct from a well-filled
tensor field $R^{(N-1)n}$ of degree $(N-1)n$ with $n\geq 1$
fulfilling $d^k R^{(N-1)n}=0$, a tensor field $K^{(N-1)(n-1)+k-1}$ (of degree $(N-1)(n-1)+k-1$) such that
$R^{(N-1)n} = d^{N-k} K^{(N-1)(n-1)+k-1}$. The construction works only for
well-filled tensors; for tensors of a different Young symmetry
type, the tensor $K$ obtained through the homotopy in the
given $N$-complex does
not fulfill $d^{N-k}K = R$, (for $d^kR=0$). 

The details will be given in \cite{D-VH}.
We shall merely display here two explicit homotopy formulas
which reveal the main points and which deals with
cohomologies effectively investigated in
the literature previously (see next section).  Consider first
in $\Omega_4(\mathbb R^D)$
a tensor $T$ in degree $3$ which is annihilated by $d^3$.  
In components,
$$\partial_{[\alpha_1} \partial_{[\beta_1} \partial_{[\gamma_1}
T_{\alpha_2]\beta_2]\gamma_2]} = 0$$
where the antisymmetries are on the $\alpha$'s, the $\beta$'s
and the $\gamma$'s.  A straightforward calculation shows that
$d^3T=0$ implies $T = d\xi$ ($\leftrightarrow
T_{\alpha \beta \gamma} = \partial_{(\alpha} \xi_{\beta \gamma)}$), 
with $\xi_{\alpha\beta}$ given by the homotopy formula
\begin{eqnarray}
\xi_{\alpha \beta}(x) &=& \int_0^1 dt \, T_{\alpha \beta \lambda} (tx)
\, x^\lambda
\nonumber \\
&+& \frac{1}{2} \int_0^1 dt \int_0^t dt' \, 
(\partial_{[\mu} T_{\alpha]\beta \lambda} (t'x) +
\partial_{[\mu} T_{\beta] \alpha \lambda} (t'x) ) \, x^\mu \, x^\lambda 
\nonumber \\
&+& \int_0^1 dt \int_0^t dt' \int_0^{t'} dt'' \,
\partial_{[\mu} \partial_{[\rho} T_{\alpha]\beta] \lambda} (t''x)
\, x^\lambda \, x^\mu \, x^\rho. \label{eq4} 
\end{eqnarray}
Thus, $H^{3}_{(3)}(\Omega_4(\mathbb R^D))=0$.  
In the homotopy formula (\ref{eq4}), not only
does the inner contraction of $T$ with $x$ appear,
but also the double contraction of $dT$
with $x x$, as well as the triple contraction
of $d^2 T$ with $x x x$. 

The second illustrative homotopy formula shows that 
$H^{4}_{(1)}(\Omega_3(\mathbb R^D))=0$.
If the tensor
$R_{\alpha_1 \alpha_2 \beta_1 \beta_2}$ of degree 4 with the symmetry
$ \begin{tabular}{|c|c|}
\hline
$\alpha_1$ & $\beta_1$ \\
\hline
$\alpha_2$ & $\beta_2$ \\
\hline
\end{tabular}\
$, (i.e. the symmetry of Riemann curvature tensor) 
is annihilated by $d$, $\partial_{[\alpha_3} 
R_{\alpha_1 \alpha_2] \beta_1 \beta_2} = 0$,
then one has $R_{\alpha_1 \alpha_2 \beta_1 \beta_2}
= 4 \partial_{[\beta_2}  \partial_{[\alpha_2} h_{\alpha_1]  \beta_1]}$
with
\begin{equation}
 h_{\alpha \beta}(x) = \int_0^1 dt \int_0^{t} dt'\  t'\  x^\lambda x^\mu R_{\alpha \lambda \beta \mu} (t'x).
\label{eq5} 
\end{equation}
Both homotopy formulas given here follow the general pattern
described in \cite{D-V2}.  The general case for arbitrary $N$, 
$n$ or $k$ in the theorem leads to homotopy formulas with the same
structure.

One can alternatively prove the theorem
by repeated use of the standard Poincar\'e lemma for differential forms,
but this appears to be more laborious for big $N$. 

Although there is no cohomology for well-filled tensors,
the cohomology is non trivial at the other tensorial degrees.
One easily verifies that the cohomology for tensors
corresponding to a single (unfilled) row is finite-dimensional
and related to the Killing tensors of Minkowski space.
The cohomology in the other cases, however, is 
generically infinite-dimensional.
One may remove it by embedding the $N$-complex
$(\Omega_N(\mathbb R^D), d)$ in a bigger $N$-complex, containing 
different symmetry types (and thus reducible tensors) in
each tensorial degree, but this will
not be done here.  Again, the details will be given in
\cite{D-VH}.

\section{Higher spin gauge fields}

The $N$-complexes (and their generalized cohomologies)
defined in this letter naturally arise
in the description of higher integer spin gauge fields.

Classical spin $S$ gauge fields (with $S\in \mathbb N$) are described by symmetric tensor
fields $h_{\alpha_1 \dots \alpha_S}$
of order $S$ and gauge transformations of the form
\begin{equation}
\delta_{\epsilon} h_{\alpha_1 \dots \alpha_S} =
\partial_{(\alpha_1} \epsilon_{\alpha_2 \dots \alpha_S)}
\label{eq6}
\end{equation}
where $\epsilon_{\alpha_2 \dots \alpha_S}$ is a symmetric tensor
of order $S-1$.
\footnote{For $S \geq 3$, the gauge parameter
is subject to the trace condition $\epsilon^{\alpha_2}_{\;
\; \; \alpha_2 \dots \alpha_{S_1}}
=0$ and for $S \geq 4$, the gauge field is subject to
the double-trace condition $h^{\alpha_1 \alpha_2}_{\; \;
\; \; \; \; \; \alpha_1 \alpha_2
\dots \alpha_{S-2}} =0$ \cite{F,SH}.  However, as observed in 
\cite{dWF,DD}, it is already
of interest to investigate the gauge symmetries without imposing
the trace conditions.}.  The curvatures
$R_{\alpha_1 \dots \alpha_S \beta_1 \dots \beta_S}$
invariant under (\ref{eq6}) contain $S$ derivatives
of the fields \cite{dWF} and are obtained from
$\partial_{\alpha_1 \dots \alpha_S} h_{\beta_1 \dots \beta_S}$
by symmetrizing according to the Young tableau with $S$ columns
and $2$ rows.

It is clear from the above definitions that $R=d^{S} h$ where $d$
is the derivative operator of the complex 
$(\Omega_{S+1}(\mathbb R^D), d)$.  Gauge invariance of the
curvature follows from $d^{S+1} = 0$.

The generalized Poincar\'e lemma (Theorem 1) implies
$H^{S}_{(S)}(\Omega_{S+1}(\mathbb R^D))=0$ which ensures that gauge
fields with zero curvatures are pure gauge.  This was directly proved in \cite{DD} for the case $S=3$.
The condition $d^{S+1}=0$ also
ensures that curvatures of gauge potentials satisfy a
generalized Bianchi
identity of the form $dR = 0$.
The generalized Poincar\'e lemma also implies 
$H^{2S}_{(1)}(\Omega_{S+1}(\mathbb R^D))=0$
which means that
conversely the Bianchi identity characterizes
the elements of $\Omega^{2S}(\mathbb R^D)$ which are curvatures of
gauge potentials. This claim for $S=2$
is the main statement of \cite{Gas}.

\section{Duality}

{}Finally, there is a generalization of Hodge duality for
$\Omega_N(\mathbb R^D)$, which is obtained by contractions of
the columns with the Kroneker tensor
$\varepsilon^{\mu_1\dots\mu_D}$ of $\mathbb R^D$. A detailed description of this duality will appear in \cite{D-VH}. 
When combined with Theorem 1, this
duality leads to another kind of results. A typical result of
this kind is the following one. Let $T^{\mu\nu}$ be a symmetric
contravariant tensor field of degree 2 on $\mathbb R^D$
satisfying $\partial_\mu T^{\mu\nu}=0$, (like e.g. the stress
energy tensor), then there is a contravariant tensor field
$R^{\lambda\mu\rho\nu}$ of degree 4 with the symmetry
$ \begin{tabular}{|c|c|}
\hline
$\lambda$ & $\rho$ \\
\hline
$\mu$ & $\nu$ \\
\hline
\end{tabular}\ 
$, (i.e. the symmetry of Riemann curvature tensor), such
that
\begin{equation}
T^{\mu\nu}=\partial_\lambda\partial_\rho R^{\lambda\mu \rho\nu}
\label{eq7}
\end{equation}
In order to connect this result with Theorem 1, define
$\tau_{\mu_1\dots\mu_{D-1}
\nu_1\dots\nu_{D-1}}=\linebreak[4]T^{\mu\nu}\varepsilon_{\mu\mu_1\dots\mu_{D-1}}
\varepsilon_{\nu\nu_1\dots \nu_{D-1}}$. Then one has
$\tau\in\Omega^{2(D-1)}_3(\mathbb R^D)$ and conversely, any
$\tau\in \Omega^{2(D-1)}_3(\mathbb R^D)$ can be expressed in
this form in terms of a symmetric contravariant 2-tensor. It is
easy to verify that $d\tau=0$ (in $\Omega_3(\mathbb R^D))$ is
equivalent to $\partial_\mu T^{\mu\nu}=0$. On the other hand,
Theorem 1 implies that $H^{2(D-1)}_{(1)}(\Omega_3(\mathbb
R^D))=0$ and therefore $\partial_\mu T^{\mu\nu}=0$ implies that
there is a $\rho\in \Omega^{2(D-2)}_3(\mathbb R^D)$ such that
$\tau=d^2\rho$. The latter is equivalent to (\ref{eq7}) with
$R^{\mu_1\mu_2\ \nu_1\nu_2}$ proportional to
$\varepsilon^{\mu_1\mu_2\dots\mu_D}
\varepsilon^{\nu_1\nu_2\dots
\nu_D}\rho_{\mu_3\dots\mu_D\nu_3\dots\nu_D}$ and one verifies
that, so defined, $R$ has the correct symmetry. That symmetric
tensor fields identically
fulfilling $\partial_\mu T^{\mu\nu}=0$ can
be rewritten as in Eq. (\ref{eq7}) has been used in
\cite{W} in the investigation of the consistent deformations of
the free spin two gauge field action. 

\section{The differential calculus for a manifold}

If the space $\mathbb R^D$ is replaced by an arbitrary $D$-dimensional smooth manifold $V$, then smooth covariant tensor fields of given Young symmetry type are still well defined and therefore the graded space $\Omega_Y(V)=\displaystyle{\oplusinf _p} \Omega^p_Y(V)$ is well defined for a sequence $Y=(Y_p)_{p\in \mathbb N}$ of Young diagrams such that $Y_p$ has $p$ cells $\forall p\in \mathbb N$ as in Section 2. In fact $\Omega_Y(V)$ is a graded module over  the algebra $C^\infty(V)$ of smooth functions and (\ref{eq3}) still defines a $C^\infty(V)$-bilinear graded product on $\Omega_Y(V)$. However now the operator $T\mapsto \partial T$ of Section 2 does not make sense; in order to give a substitute for it, one must choose a linear connection on $V$ and replace $\partial$ by the corresponding covariant derivative $\nabla$. One then generalizes $d$ by $d_\nabla=\mathbf Y\circ \nabla$, i.e. formula (2) by
\begin{equation}
d_\nabla=\mathbf Y_{p+1}\circ \nabla:\Omega^p_Y(V)\rightarrow \Omega^{p+1}_Y(V)
\label{eq8}
\end{equation}
which defines again a first order differential operator on $\Omega_Y(V)$. This operator $d_\nabla$ is again homogeneous of degree 1 but now, due to the torsion and the curvature of $\nabla$, Lemma 1 is not true.  In fact Lemma 1 merely applies at the level of symbols; more precisely one has the following: {\sl Let $N$ and $Y$ satisfy the assumptions of Lemma 1, then $(d_\nabla)^N$ is a differential operator of order smaller or equal to $N-1$ and, if furthermore $\nabla$ is torsion-free, then the order of $(d_\nabla)^N$ is smaller or equal to $N-2$.} In the case $N=2$, if $\nabla$ is torsion free, $(d_\nabla)^2=0$ follows from the first Bianchi identity; however in this case $d_\nabla$ coincides, as well known, with the exterior differential $d$ which is well defined in local coordinates by (\ref{eq1}).

\section*{Acknowledgements}
M.H is grateful to the ``Laboratoire de Physique Th\'eorique
de l'Universit\'e Paris XI" for kind hospitality while
this work was carried out.  His research was partly supported
by the ``Actions de
Recherche Concert{\'e}es" of the ``Direction de la Recherche
Scientifique - Communaut{\'e} Fran{\c c}aise de Belgique", by
IISN - Belgium (convention 4.4505.86) and by
Proyectos FONDECYT 1970151 and 7960001 (Chile).

\end{document}